\documentclass[12pt]{amsart}
\usepackage{amscd, amssymb, epic, eepic}

\theoremstyle{plain}

\newtheorem{lemma}{Lemma}

\newtheorem{theorem}{Theorem}
\newtheorem*{main}{Main Theorem}
\newtheorem{corollary}{Corollary}
\newtheorem*{claim*}{Claim}
\newtheorem{conjecture}{Conjecture}
\newtheorem*{question*}{Question}
\newtheorem{remark}{Remark}

\theoremstyle{definition}
\newtheorem{definition}{Definition}
\newtheorem*{definition*}{Definition}
\newtheorem{convention}{Convention}
\newtheorem*{convention*}{Convention}

\theoremstyle{remark}
\newtheorem*{remark*}{Remark}
\newtheorem*{remarks*}{Remarks}

\newtheorem*{example*}{Example}
\newtheorem*{acknowledgement}{Acknowledgement}

\newcommand{\PP}{\mathbb{P}}

\newcommand{\QQ}{\mathbb{Q}}

\newcommand{\Mbar}{{\overline{M}}}

\newcommand{\p}{\partial}

\newcommand{\ft}{\operatorname{ft}}
\newcommand{\st}{\operatorname{st}}

\newcommand{\la}{\langle}
\newcommand{\ra}{\rangle}

\title[Invariance of tautological equations I]
{Invariance of tautological equations I: conjectures and applications}

\author{Y.-P.~Lee}

\address{Department of Mathematics \\
        University of Utah \\
        Salt Lake City, Utah 84112-0090\\
        U.S.A.}
\email{yplee@math.utah.edu}

\thanks{Research partially supported by NSF and AMS Centennial Fellowship}

\begin{document}

\begin{abstract}
The main goal of this paper is to introduce a set of conjectures on the
relations in the tautological rings.
In particular, this framework gives an efficient algorithm to calculate
all tautological equations using only finite dimensional linear algebra.
Other applications include
the proofs of \emph{Witten's conjecture} on the relations
between higher spin curves and Gelfand--Dickey hierarchy and
\emph{Virasoro conjecture} for target manifolds with conformal semisimple
quantum cohomology, both for genus up to two.
\end{abstract}

\maketitle

\section{Introduction}

\subsection{The tautological rings of moduli spaces of curves}
Two basic references for tautological rings are \cite{HM} and \cite{rV}.

Let $\Mbar_{g,n}$ be the moduli stacks of stable curves.
$\Mbar_{g,n}$ are proper, irreducible, smooth Deligne--Mumford stacks.
The Chow rings $A^*(\Mbar_{g,n})$ over $\QQ$ are isomorphic to
the Chow rings of their coarse moduli spaces.
The tautological rings $R^*(\Mbar_{g,n})$ are subrings of
$A^*(\Mbar_{g,n})$, or subrings of $H^{2*}(\Mbar_{g,n})$ via cycle maps,
generated by some ``geometric classes'' which will be described below.

\begin{convention}
All Chow/Cohomology/tautological rings are over $\QQ$.
\end{convention}

The first type of geometric classes are the \emph{boundary strata}.
$\Mbar_{g,n}$ have natural stratification by topological types.
The strata can be conveniently presented by their (dual) graphs, which
can be described as follows.
To each stable curve $C$ with marked points, one can associate a dual graph
$\Gamma$.
Vertices of $\Gamma$ correspond to irreducible components.
They are labeled by their geometric genus.
Assign an edge joining two vertices each time the two components intersect.
To each marked point, one draws an half-edge incident to the vertex,
with the same label as the point.
Now, the stratum corresponding to $\Gamma$ is the closure of the subset of all
stable curves in $\Mbar_{g,n}$ which have the same topological type as $C$.

The second type of geometric classes are the Chern classes of
tautological vector bundles. These includes cotangent classes $\psi_i$,
Hodge classes $\lambda_k$ and $\kappa$-classes $\kappa_l$.
(See \cite{HM}, \cite{rP}.)

To give a precise definition of the tautological rings,
some natural morphisms between moduli stacks of curves will be used.
The \emph{forgetful morphisms}
\begin{equation} \label{e:ft}
 \ft_i: \Mbar_{g,n+1} \to \Mbar_{g,n}
\end{equation}
forget one of the $n+1$ marked points.
The \emph{gluing morphisms}
\begin{equation} \label{e:gl}
 \Mbar_{g_1,n_1 +1} \times \Mbar_{g_2, n_2 +1} \to \Mbar_{g_1 + g_2, n_1+n_2},
 \quad \Mbar_{g-1,n+2} \to \Mbar_{g,n},
\end{equation}
glue two marked points to form a curve with a new node.
Note that the boundary strata are the images (of the repeated applications)
of the gluing morphisms, up to factors in $\QQ$ due to automorphisms.

\begin{definition}
The system of tautological rings $\{R^*(\Mbar_{g,n})\}_{g,n}$ is the
smallest system of $\QQ$-unital subalgebra (containing classes of
type one and two, and is) closed under the forgetful and gluing morphisms.
\end{definition}

\begin{remarks*}
(i) The phrase ``(containing classes of type one and two, and is)'' can
be removed from the definition as the type one classes can be generated
by the fundamental classes of the boundary strata (units), and the type two
classes can be generated by the type one classes under the natural morphisms.

(ii) The set of Hodge classes can be reconstructed \emph{linearly} from
the set of $\psi$-classes, $\kappa$-classes and boundary strata. 
Therefore, the Hodge classes
can be omitted in the discussion of tautological rings.
\end{remarks*}

Since the tautological rings $\{R^*(\Mbar_{g,n})\}_{g,n}$ are defined by
generators, it leaves us the task of finding the \emph{relations} in order to
understand their structures. Mumford \cite{dM}, Getzler \cite{eG1}
and Faber \cite{cF} are among the pioneers in this direction.

\subsection{The conjectures} \label{s:1.2}

For each dual graph $\Gamma$, one can decorate the graph by assigning
a monomial, or more generally a polynomial,
of $\psi$ to each half-edge and $\kappa$ classes to each vertex.
The tautological classes in $R^k(\Mbar_{g,n})$ can be
represented by $\QQ$-linear combinations of \emph{decorated graphs},
which will be called simply graphs by abusing the notations.

\begin{convention}
The graphical notations might be a little different from some authors'.
\emph{The edges here are considered as ``gluing'' two half-edges.}
Therefore, the consideration of automorphisms will lead to discrepancy of
a constant factor.
\end{convention}

Define the operations $\mathfrak{r}_l$ on the spaces of decorated graphs
$\{ \Gamma \}$.

\begin{itemize}
\item \textbf{Cutting edges.} Cut one edge and create two new half-edges.
Label two new half-edges with $i,j \notin \{1,2,\ldots,n\} $ 
in two different ways.
Produce a formal sum of 4 graphs by decorating extra $\psi^l$ to $i$-labeled 
new half-edges with coefficient $1/2$ and by decorating extra $\psi^l$ to 
$j$-labeled new half-edges with coefficient $(-1)^l/2$  . 
(By ``extra'' decoration we mean that $\psi^l$ is multiplied by
whatever decorations which are already there.)
Produce more graphs by proceeding to the next edge.
Retain only the stable graphs. 
Take formal sum of these final graphs.

\item \textbf{Genus reduction.}
For each vertex, produce $l$ graphs. Reduce the genus of this given
vertex by one, add two new half-edges.
Label two new half-edges with $i,j$ and decorate them
by $\psi^{l-1-m}, \psi^m$ (respectively) where $0 \le m \le l-1$.
Do this to all vertices, and retain only the stable graphs.
Take formal sum of these graphs with coefficient $\frac{1}{2} (-1)^{m+1}$.

\item \textbf{Splitting vertices.} Split one vertex into two.
Add one new half-edge to each of the two new vertices.
Label them with $i,j$ and decorate them by $\psi^{l-1-m}, \psi^m$
(respectively) where $0 \le m \le l-1$.
Produce new graphs by splitting the genus $g$ between the two new vertices
($g_1, g_2$ such that $g_1+g_2=g$),
and distributing to the two new vertices the (old) half-edges which
belongs to the original chosen vertex, in all possible ways.
The $\kappa$-classes on the given vertex are split between the two 
new vertices in a way similar to the half-edges.
That is, consider each monomial of the $\kappa$-classes 
$\kappa_{l_1} \ldots \kappa_{l_p}$ on the split vertex as labeled by 
$p$ special half-edges. 
When the vertex splits, distribute the $p$ special edges in all possible way.
Do this to all vertices, and retain only the stable graphs.
Take formal sum of these graphs with coefficient $\frac{1}{2} (-1)^{m+1}$.

\end{itemize}

\begin{remarks*}
(i) It is not very difficult to see that the two new half-edges are
symmetric for $l$ odd and anti-symmetric for $l$ even.


(ii) The (output) graphs might be disconnected. The stable graphs here mean
that each connected component is stable (and of non-negative dimension).

(iii) These operations in terms of another (equivalent) notation 
(called \emph{gwi}'s) will be given in Section~3.1.

(iv) The $\lambda$-classes do not enter the discussion as they can be 
reconstructed from other classes \cite{dM} \cite{FP}. 
To include them explicitly, 
one may apply the elementary splittings of the Hodge bundle
\[
  0 \to \mathbb{E}_{g-1} \to \iota^* \mathbb{E}_g \to \mathcal{O}_{\Delta_0} \to 0
\]
to the genus-reduced vertices and
\[
  0 \to \mathbb{E}_{g'} \to \iota^* \mathbb{E}_g \to \mathbb{E}_{g-g'} \to 0
\]
to the split vertices.
\end{remarks*}

\begin{definition}
Define operators $\mathfrak{r}_l(\Gamma)$ from the vector spaces of
decorated graphs of codimension $k$ on $\overline{M}_{g,n}$ to
those of codimension $k+l-1$ on $\overline{M}^{\bullet}_{g-1,n+2}$
to be the final sum of the decorated graphs (with $\QQ$-coefficients)
produced from the above three operations. 
Here $\bullet$ stands for possibly disconnected curves.
Note that the arithmetic genus for disconnected curve is defined to be
\[
  g(C) := \sum_{i=1}^d g(C_i) - d +1,
\]
where $C_i$ are connected components of $C$, $C= \amalg_{i=1}^d C_i$.
\end{definition}

Let
\[
  \sum_i c_i \Gamma_i =0
\]
be a tautological equation in codimension $k$ classes in $\overline{M}_{g,n}$.
The triple $(g,n,k)$ will be used.
\begin{conjecture}
\footnote{This conjecture will be proved in \cite{ypL2}.}
\begin{equation} \label{e:inv}
 \mathfrak{r}_l (\sum_i c_i \Gamma_i) =0
\end{equation}
for all $l$, modulo the tautological equations for $(g', n')$ for which
\begin{itemize}
\item $g' < g$, or
\item $g'=g$ and $n'<n$.
\end{itemize}
\end{conjecture}

The set of equations \eqref{e:inv} will be called the 
\emph{$R$-invariance} equations.
The reason of this name comes from its relation to Gromov--Witten theory
and will be explained in \cite{ypL2}.
See also \cite{ypL1} for a brief account.

\begin{conjecture}
Let $E=\sum_i c_i \Gamma_i$ be a given linear combination of codimension $k$
tautological classes in $\Mbar_{g,n}$ and $k < 3g-3+n$.

If $\mathfrak{r}_l (E) =0$ for all $l$, modulo tautological equations in
$\overline{M}_{g',n'}$ for $(g', n')$ satisfying the above inductive
conditions, then $E=0$ is a tautological equation.
\end{conjecture}

\begin{remark}
(i) By a theorem of Graber and Vakil \cite{GV1}, when $k=3g-3+n$,
i.e.~top codimension,
\[
 R^{top}(\Mbar_{g,n})=\QQ.
\]
Therefore $R^{top}(\Mbar_{g,n})$ are considered well understood and
used as part of inductive data.

(ii) Conjecture~1 means that the linear operators at the level of decorated 
graphs actually descend to operators at the level of tautological classes.
That is,
\begin{equation} \label{e:d}
 \mathfrak{r}_l: R^{k}(\Mbar^{\bullet}_{g,n}) \to 
R^{k+l-1}(\Mbar^{\bullet}_{g-1,n+2})
\end{equation}
is well-defined.
Conjecture~2 further asserts that the combination of these operators 
is injective.
\end{remark}

\begin{conjecture}
Conjecture~2 will produce \emph{all} tautological equations inductively.
\end{conjecture}

\subsection{Results in low genus} \label{s:1.3}

\begin{main} \cite{GL, AL1, AL2}
The Invariance Conjectures hold for genus zero, genus one, and
$(g,n,k)=(2,1,2), (2,2,2), (2,3,2), (3,1,3)$.
\end{main}

Genus zero case: $R^k(\Mbar_{0,n}) = A^k(\Mbar_{0,n})$ essentially
due to S.~Keel \cite{sK}.
There are two sets of tautological equation in genus zero,
namely the genus zero topological recursion relations (TRR) and
WDVV equations. WDVV equations are consequences of TRRs.
TRR is based upon the fact that $\psi_i=0$ on $\Mbar_{0,3}$.
WDVV equation is based upon the fact that $A^1(\Mbar_{0,4})=\QQ$.


The non-trivial known tautological equations in $g \le 3$ are
\begin{enumerate}
\item $(g,n,k)=(1,4,2)$ Getzler's genus one equation (\cite{eG1}).
\item $(g,n,k)=(2,1,2)$ Mumford-Getzler's equation;
\item $(g,n,k)= (2,2,2)$ Getzler's genus two equation (\cite{eG2}).
\item $(g,n,k)= (2,3,2)$ Belorousski--Pandharipande's equation (\cite{BP}).
\item $(g,n,k)=(3,1,3)$ A new relation (\cite{AL2} \cite{KL}).
\end{enumerate}

Genus one case will be discussed in Section~\ref{s:3}.

Genus two and three.
The checking of cases (2-4) is carried out in \cite{AL1}.
The calculation of (5) via Conjecture~2 is done in \cite{AL2}.
An equivalent form of the $(3,1,3)$ equation was independently discovered by
Kimura--Liu \cite{KL} using a completely different method.
(Equivalent in the sense that they use different vectors in $R^3(\Mbar_{3,1})$.)
By Getzler's Betty number and Hodge polynomial calculations
in \cite{eG2} and \cite{eGL}, these equations are the only tautological
equations in the above cases.

\begin{remark} \label{r:1}
(i) Here the term ``the only tautological equations'' should be taken with
a grain of salt. The forgetful and gluing morphisms \eqref{e:ft} \eqref{e:gl}
induce relations in $R^k(\Mbar_{g,n})$ from lower (inductive) classes.
For example, Getzler's equation in $(g=1,n=4,k=2)$ induces an equation
in $(2,2,3)$ by gluing two of the marked points;
any equations in $\Mbar_{g,n}$ will induce equations in $\Mbar_{g,n+m}$.
These induced equations will be taken into inductive data as well.
The goal, of course, is to find new equations.

(ii) It is easy to see that if one equation is $R$-invariant, all
induced equations are.
\end{remark}

\subsection{Relations to Gromov--Witten theory}
This will be the topic for a forthcoming paper \cite{ypL2},
so the discussion will be brief and necessarily not precise.

\subsubsection{Motivation of the conjectures}
The motivation of the above conjectures come from a study of Givental's
axiomatic Gromov--Witten theory \cite{aG1} \cite{aG2} \cite{aG4}.
Givental has discovered some remarkable structures of the ``moduli spaces'' of
the semisimple axiomatic Gromov--Witten theories (or Frobenius manifolds).
The ``moduli space'' of a given rank (i.e.~the dimension of the Frobenius
manifolds) is a ``homogenius space'' of a ``quantized loop group''.

(i) The tautological equations hold for all geometric Gromov--Witten theories
due to the fact that there is a natural stabilization morphism
\[
 \st: \Mbar_{g,n}(X,\beta) \to \Mbar_{g,n}
\]
from the moduli spaces of stable maps to moduli spaces of stable curves.
It is natural to expect that the tautological equations hold for all
axiomatic Gromov--Witten theories. Therefore, the tautological equations
should be ``invariant'' under the action of the ``quantized loop groups''.
This translates eventually to Conjecture~1.

(ii) Conversely, it is also expected that all ``universal equations'' of
(axiomatic) Gromov--Witten invariants should be induced from
the tautological equations on moduli spaces of curves.
The term ``universal equations'' means that the equations are valid
and of the same form for all theories.
In other words, they are invariant under the action of the ``quantized
loop groups''.
This gives Conjectures~2 and 3.

\subsubsection{Applications in axiomatic Gromov--Witten theory}
Combining the results of this paper with two separate papers \cite{ypL2}
\cite{AL1},
the last one joint with Arcara, the following results are proved.

\begin{theorem}
(i) All tautological equations discussed in Section~\ref{s:1.3}
hold for Givental's axiomatic Gromov--Witten invariants.

(ii) The Virasoro Conjecture for semisimple conformal Frobenius manifolds
holds up to genus two.

(iii) The Witten's conjecture on the higher spin curves and Gelfand--Dickey
hierarchies holds up to genus two.
\end{theorem}

\begin{remarks*}
(i) In a joint work with Givental \cite{GL}, it is shown that the conformal
condition in Theorem (ii) can be removed for $g=1$.

(ii) Theorem (ii) was independently proved by X.~Liu \cite{xL}.

(iii) The above Theorems are built upon many other authors' results, and will
be discussed in \cite{ypL2}.
\end{remarks*}

\begin{acknowledgement}
Special thanks are due to my collaborators D.~Arcara (\cite{AL1,AL2}) 
and A.~Givental (\cite{GL}) for many useful discussions 
during our joint projects related to this paper;
to R.~Cavalieri and R.~Vakil for their interest and clarification to the 
conjectures;
to participants of WAGS Spring 2004 and 
Moduli Space of Curves and Gromov-Witten Theory Workshop
for suggestions in the expositions.
\end{acknowledgement}

\section{Algorithm of finding tautological equations}

\subsection{Finiteness}

Assuming Invariance Conjectures, equation \eqref{e:d} says that
\[
 \mathfrak{r}_l: R^{k}(\Mbar_{g,n}) \to 
R^{k+l-1}(\Mbar^{\bullet}_{g-1,n+2})
\]
is injective.
In fact, the number of connected components in the image 
can go up by at most one.
Therefore, if the resulting decorated graphs are connected, the genus
must be reduced by one.
If disconnected, then either genus or the number of marked points
(external half-edges) is reduced.

\begin{lemma} \label{l:1}
$\mathfrak{r}_l$ reduces the dimension of $\Gamma$ by $l$.
\end{lemma}

\begin{proof}
The proof is straightforward case-by-case study.
In the first step, the dimension remains the same after cutting an edge.
Decorating with extra $\psi^l$ reduces the dimension by $l$.
In the second step, reducing the genus and adding two half-edges change
the dimension by $-3+2=1$. Decorating with $\psi^m$ and $\psi^{l-1}$ reduces
the dimension by $l-1$.
In the third step, the vertex splitting reduces the dimension by $3$.
Adding two half-edges increases dimension by $2$. Decoration reduces
the dimension by $l-1$.
\end{proof}

\begin{remark} \label{r:2}
(i) By Lemma~\ref{l:1},
$\mathfrak{r}_l \Gamma =0$ when $k+l > \dim \overline{M}_{g,n}$.
Therefore, one only has to check finite number of $l$'s.

(ii) There are only finitely many $(g',n',k')$ involved in checking the
validity of invariance equation \eqref{e:inv} for a given $(g,n,k)$ and $l$.
This observation can be easily verified from the definition
of $\mathfrak{r}_l$.
Therefore, it is a finite calculation to check \eqref{e:inv} for any quadruple
$(g,n,k,l)$.
\end{remark}

From the above discussions, one knows that the

\begin{corollary}
There is an efficient algorithm to check Conjectures~1 and 2 for any given
$(g,n,k)$.
\end{corollary}

\subsection{The algorithm}
Invariance Conjectures have been used to ``re-discover'' all known 
tautological equations and discover a new one.
Our basic strategy of finding the tautological equations for a given
$(g,n,k)$ is the following. Assuming all tautological equations are known
for $(g',n',k')$ satisfying the inductive conditions in Conjecture~1.

\begin{enumerate}
\item Find all tautological classes $\Gamma_i$ in $R^k(\Mbar_{g,n})$.
Remove linear dependent classes from the induced equations as explained in
Remark~\ref{r:1}.
Let $\{ \Gamma_i \}_{i \in I}$ be a set of remaining vectors in the
$\QQ$-linear space $R^k(\Mbar_{g,n})$.
Let
\begin{equation*}
 E = \sum_{i \in I} c_i \Gamma_i
\end{equation*}
be a general element in
$R^k(\Mbar_{g,n})$ with \emph{unknown} coefficients $c_i$.

\item Apply the $R$-invariance condition \eqref{e:inv}
$\mathfrak{r}_l (E)=0$, for $l=1,\ldots, 3g-3+n-k$.
For each $l$, $\mathfrak{r}_l(\Gamma)$ will be a $\QQ$-linear combination
of (disjoint unions of) $\Gamma'_j$ as classes in $R^{k'}(\Mbar_{g',n'})$,
which are known by assumption. In particular, one knows all the relations
between $\Gamma'_j$'s.

\item In each $(g',n',k')$ pick a \emph{basis} of $R^{k'}(\Mbar_{g',n'})$.
The part of the output of $\mathfrak{r}_l (E)=0$ in $(g',n',k')$ implies each
component of the basis vanishes individually.
Since all operations involved are \emph{linear}, the vanishing gives
linear equations in $c_i$.

\item The above step produces enough linear equations on $c_i$'s to determine
them completely up to a few free variables, say $c_1$ and $c_2$.
Write all other $c_j$'s in terms of $c_1$ and $c_2$.
$E=\sum_{i \in I} c_i \Gamma_i $  becomes
\[
 E= c_1 ( \sum_{j_1 \in J_1} d_{j_1} \Gamma_{j_1})
 + c_2 (\sum_{j_2 \in J_2} d_{j_2} \Gamma_{j_2}).
\]
The output equations are
\[
  \sum_{j_1 \in J_1} d_{j_1} \Gamma_{j_1} =0, \quad \textrm{and} \quad
  \sum_{j_2 \in J_2} d_{j_2} \Gamma_{j_2} =0,
\]
where $d_j$ are output constants.
Note that $J_1$ and $J_2$ are not necessarily disjoint

\item If all $c_i$ have to vanish after the above steps, then there is no
(new) tautological equations for $(g,n,k)$.
\end{enumerate}

\begin{remarks*}
(i) The Step (4) will necessarily leave at least one free variable, as
any equation $E=0$ holds after a constant multiple.

(ii) In fact, the output of the $R$-invariance condition always highly
\emph{over-determines} the unknown coefficients $c_i$.
To be able to solve $c_i$ is usually a sign of correct calculations.

(iii) In all cases computed, $l=1$ is enough. That is, $\mathfrak{r}_1 (E)=0$
already generate enough linear equations among $c_i$ to determine them
completely. It is not known whether this will hold in general.
In particular, one might ask whether $\mathfrak{r}_1 (E)=0$ implies
$\mathfrak{r}_l (E)=0$ for $l \ge 2$.
We don't have intuition leading to a guess.
\end{remarks*}

As stated above, the above algorithm has given a \emph{uniform} method of
deriving all tautological equations, which were originally derived using many
different methods.
Furthermore, it is theoretically possible to program this algorithm in order 
to discover more tautological equations via robots.
The only ingredient in this algorithm is linear algebra
(and the efforts to obtain induced equations as explained in Remark~\ref{r:1}).
However, the dimension of $R^k(\Mbar_{g,n})$ and
the number of elements in the set $\{ \Gamma_i \}_{i \in I}$ can grow.
So the finite dimensional linear algebra problem in question is not trivial.

\begin{remark}
An alternative way to the above (more satisfactory) algorithm goes through
the following procedure.
\begin{itemize}
\item Calculate the rank of $R^k(\Mbar_{g,n})$ to 
see if there is any new equation.
\item If there is one, then apply invariance condition equation \eqref{e:inv} 
to obtain the coefficients of the equation.
\end{itemize}
This is actually the way which is employed to prove the derived tautological
equations. However the first step is usually not elementary.
\end{remark}

\section{Proof of Conjectures for $g=1$} \label{s:3}

\subsection{Notations}
In this subsection, some notations will be introduced for future reference.

The following notations will be used to denote the same object.
\begin{enumerate}
\item Tautological classes.
\item Generic curves and chern classes.
\item Decorated graphs.
\item Gwi's (defined below).
\end{enumerate}

(1) $\Leftrightarrow$ (2).
For each tautological class, one may draw a generic curve of the given
topological type, label marked points, and decorate the marked points and
both sides of nodes with monomials of $\psi$ classes and the 
components by $\kappa$ classes.

(1) $\Leftrightarrow$ (3). This has been explained in Section~1.1-1.2.

(3) $\Leftrightarrow$ (4).
Given a decorated graph $\Gamma$.
\begin{itemize}
\item For the vertices of $\Gamma$ of genus $g_1, g_2, \ldots$,
assign a product of ``brackets''
$\langle \rangle_{g_1} \langle \rangle_{g_2} \ldots$.
\item Assign each half-edge a symbol $\p^{*}$.
The external half-edges use numeric super-indices $\p^1,\p^2,\ldots, \p^n$,
corresponding to their labeling.
The two new half-edges use $\p^i,\p^j$.
For each pair of half-edges coming from one and the same edge, the same 
super-index $\mu_i$ will by used.
Otherwise, all half-edges should use different super-indices.
\item For each a given vertex $\langle \rangle_{g}$ with $m$ half-edges,
$n$ external half-edges, and say 2 new half-edges,
an insertion is placed in the vertex
$\langle \p^i \p^j \prod_{a=1}^n \p^a \prod_{b=1}^m \p^{\mu_b} \rangle_{g}$.
\item For each decoration to a half-edge by chern classes $c$,
assign a subindex to the corresponding half-edge $\p^{\mu}_c$.
\end{itemize}
The output is called a \emph{gwi}.

\begin{remarks*}
(i) Gwi is so named for its relations to Gromov--Witten invariants.
Note however that gwi here stands for a vector in a $\QQ$-algebra, rather
than a rational number.
The actual relation is rather involved and will be discussed
in details in \cite{ypL2}.
Gwi notations is used in the proof of this paper, mainly for typesetting
convenience.

(ii) The Convention on graphs adopted in Section~\ref{s:1.2} is meant to match
the graphical and gwi notations.
They often differ by a constant in some authors' conventions.
(See, e.g.~\cite{eG1}.)
\end{remarks*}

\begin{example*}
Let $\Gamma$ be the following graph. 

\vspace{20pt}

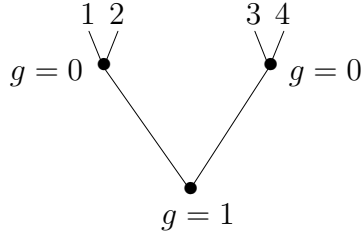
\begin{center}

\setlength{\unitlength}{0.0005in}
\begingroup\makeatletter\ifx\SetFigFont\undefined%
\gdef\SetFigFont#1#2#3#4#5{%
  \reset@font\fontsize{#1}{#2pt}%
  \fontfamily{#3}\fontseries{#4}\fontshape{#5}%
  \selectfont}%
\fi\endgroup%
{\renewcommand{\dashlinestretch}{30}
\begin{picture}(3142,2181)(0,-10)
\path(675,1779)(1575,504)(2400,1779)
\path(525,2154)(675,1779)(825,2154)
\path(2250,2154)(2400,1779)(2550,2154)
\put(600,1730){\makebox(0,0)[lb]{\smash{{{\SetFigFont{12}{14.4}{\rmdefault}{\mddefault}{\updefault}$\bullet$}}}}}
\put(2330,1730){\makebox(0,0)[lb]{\smash{{{\SetFigFont{12}{14.4}{\rmdefault}{\mddefault}{\updefault}$\bullet$}}}}}
\put(1500,440){\makebox(0,0)[lb]{\smash{{{\SetFigFont{12}{14.4}{\rmdefault}{\mddefault}{\updefault}$\bullet$}}}}}
\put(-300,1650){\makebox(0,0)[lb]{\smash{{{\SetFigFont{12}{14.4}{\rmdefault}{\mddefault}{\updefault}$g=0$}}}}}
\put(2600,1650){\makebox(0,0)[lb]{\smash{{{\SetFigFont{12}{14.4}{\rmdefault}{\mddefault}{\updefault}$g=0$}}}}}
\put(1275,150){\makebox(0,0)[lb]{\smash{{{\SetFigFont{12}{14.4}{\rmdefault}{\mddefault}{\updefault}$g=1$}}}}}
\put(430,2229){\makebox(0,0)[lb]{\smash{{{\SetFigFont{12}{14.4}{\rmdefault}{\mddefault}{\updefault}$1$}}}}}
\put(730,2229){\makebox(0,0)[lb]{\smash{{{\SetFigFont{12}{14.4}{\rmdefault}{\mddefault}{\updefault}$2$}}}}}
\put(2155,2229){\makebox(0,0)[lb]{\smash{{{\SetFigFont{12}{14.4}{\rmdefault}{\mddefault}{\updefault}$3$}}}}}
\put(2455,2229){\makebox(0,0)[lb]{\smash{{{\SetFigFont{12}{14.4}{\rmdefault}{\mddefault}{\updefault}$4$}}}}}
\end{picture}
}

\end{center}
\vspace{10pt}

The corresponding gwi is:

\[   
   \la \p^1 \p^2 \p^{\mu} \ra  \la \p^3  \p^4  \p^{\nu} \ra
    \la \p^{\mu} \p^{\nu} \ra_1 .
\]

The \emph{cutting edges} operation for $l=1$ produces
\[
  \la \p^1 \p^2 \p^i \ra \la \p^3 \p^4 \p^{\nu} \ra \la \p^j_1 \p^{\nu} \ra_1
 +\la \p^1 \p^2 \p^{\mu}\ra \la \p^3 \p^4 \p^i \ra \la \p^{\mu} \p^j_1 \ra_1.
\]
Note that $\la \p^1 \p^2 \p^i_1 \ra$ has dimension $-1$. 
Therefore the corresponding graphs are removed from the output.
Also, the $i,j$ are symmetric. Hence, a factor of $2$ is placed in front
of a term instead of adding an additional term with $i,j$ exchanged.
\end{example*}

More generally, for any edge in an expression
$\la \ldots \p^{\mu} \ldots \p^{\mu} \ldots \ra_g$, 
where the middle $\ldots$ might contain edges like $\ldots \ra_h \la \ldots$,
the cutting edge operation for $\mathfrak{r}_l$ can be written as
\begin{equation} \label{e:ce}
 \begin{split}
 &\la \ldots \p^{\mu} \ldots \p^{\mu} \ldots \ra_g  \mapsto \\
 & \frac{1}{2} \left(\la \ldots \p^i_l \ldots \p^j \ldots \ra_g 
  + (-1)^{l-1} \la \ldots \p^i \ldots \p^j_l \ldots \ra_g \right)\\
 +&\frac{1}{2} \left(\la \ldots \p^j \ldots \p^i_l \ldots \ra_g
  + (-1)^{l-1}  \la \ldots \p^j_l \ldots \p^i \ldots \ra_g \right).
 \end{split}
\end{equation}

The \emph{genus reduction} for $\mathfrak{r}_l$  produces
\begin{equation} \label{e:gr}
  \langle \p^{i_1}_{k_1} \p^{i_2}_{k_2} \ldots \rangle_g \mapsto
  \frac{1}{2} \sum_{m=0}^{l-1} (-1)^{m+1}
   \la \p^i_{l-1-m} \p^j_m \p^{i_1}_{k_1} \p^{i_2}_{k_2} \ldots \ra_{g-1}.
\end{equation}

The \emph{splitting vertices} operation for $l=1$ produces
\begin{equation} \label{e:sv}
  \langle \p^{i_1}_{k_1} \p^{i_2}_{k_2} \ldots \rangle_g \mapsto
  \frac{1}{2} \sum_{m=0}^{l-1} (-1)^{m+1}
   \sum_{g_1+g_2=g} \p^{i_1}_{k_1} \p^{i_2}_{k_2} \ldots
   ( \langle \p^i_{l-1-m} \rangle_{g_1} \langle \p^j_m \rangle_{g_2} ).
\end{equation}
The symbol of $\p^i_k$ considered as a ``linear operator'' on graphs 
is defined by Leibniz rule
\[
 \p^i_k ( \la \p^i_{l-1-m} \ra_{g_1} \la \p^j_m \ra_{g_2} )=
  \la \p^i_k \p^i_{l-1-m} \ra_{g_1} \la \p^j_m \ra_{g_2} +
  \la \p^i_{l-1-m} \ra_{g_1} \la \p^i_k \p^j_m \ra_{g_2}.
\]

\begin{convention}
In the calculations when the $\kappa$ classes are not needed,
a simplified notation is used
\begin{equation*}
 \p^{\mu}_{k} := \p^{\mu}_{\psi^k}, \quad \p^{\mu} := \p^{\mu}_{0}.
\end{equation*}
To further simplify the notations,
\[
 \langle \ldots \rangle := \langle \ldots \rangle_{g=0}.
\]
\end{convention}

\subsection{Reduction to Getzler's equation}
By a result of E.~Getzler (unpublished), the only (new) tautological equations
in genus one are his equation in $(g=1,n=4,k=2)$ \cite{eG1} and
TRR in $(g=1,n=1,k=1)$.
$g=1$ TRR expresses $\psi$-classes as boundary divisors
\begin{equation} \label{e:4}
 \langle \p^i_1 \rangle_1
 = \frac{1}{24} \langle \p^i \p^{\mu} \p^{\mu} \rangle.
\end{equation}
Since it is an equation in top codimension, it is considered as part of the
inductive data (and satisfies invariance equation by Lemma~\ref{l:1}).
Therefore, one only has to check Conjectures~1 and 2 for Getzler's equation.

\subsection{Getzler's equation in $(g,n,k)=(1,4,2)$}
The calculation here is reproduced from a joint work with A.~Givental
\cite{GL}, by referees' demand.

\begin{theorem} \cite{GL}
Getzler's equation is the only (new) codimension $2$ equation in $\Mbar_{1,4}$
which satisfies the invariance equation \eqref{e:inv}.
Furthermore, invariance equation determines the coefficients of Getzler's
equation up to a common scaling constant.
\end{theorem}

The proof is divided into the following 6 Steps.

\subsubsection{Step 1: Enumerate all boundary strata}
There are 9 codimension two boundary strata in $\Mbar_{1,4}$,
when the ordering of the 4 external half-edges are ignored. Equivalently,
one may symmetrize the 4 external half-edges by the $S_4$ permutation.
Let $\p^1, \p^2, \p^3, \p^4$ denote the 4 external half-edges.
A general element can be written as

\[
 \begin{split}
  E = \sum_{\textrm{$S_4$ permutation}} \bigg(
   &c_1 \la \p^1 \p^2 \p^{\mu} \ra  \la \p^3  \p^4  \p^{\nu} \ra
    \la \p^{\mu} \p^{\nu} \ra_1 \\
  +&c_2 \la \p^1 \p^2 \p^{\mu} \ra \la \p^3  \p^{\mu}  \p^{\nu} \ra
    \la \p^4 \p^{\nu} \ra_1 \\
  +&c_3 \la \p^1 \p^2 \p^{\mu} \ra
   \la \p^3 \p^4  \p^{\mu} \p^{\nu} \ra  \la \p^{\nu} \ra_1 \\
  +&c_4 \la \p^1 \p^2 \p^3 \p^{\mu} \ra
   \la \p^4  \p^{\mu} \p^{\nu} \ra \la \p^{\nu} \ra_1 \\
  +&c_5 \la \p^1 \p^2 \p^3 \p^{\mu} \ra
    \la \p^4 \p^{\mu} \p^{\nu} \p^{\nu} \ra  \\
  +&c_6 \la \p^1 \p^2 \p^3 \p^4 \p^{\mu} \ra
    \la \p^{\mu} \p^{\nu} \p^{\nu} \ra \\
  +&c_7 \la \p^1 \p^2 \p^{\mu} \p^{\nu} \ra
    \la \p^3 \p^4 \p^{\mu} \p^{\nu} \ra \\
  +&c_8 \la \p^1 \p^2 \p^{\mu} \ra
    \la \p^3 \p^4 \p^{\mu} \p^{\nu}  \p^{\nu} \ra  \\
  +&c_9 \la \p^1 \p^{\mu} \p^{\nu} \ra
    \la \p^2 \p^3 \p^4 \p^{\mu} \p^{\nu} \ra \bigg). \\
 \end{split}
\]

\subsubsection{Step 2: Apply Invariance Equation}
Invariance equation produces, after genus zero TRR,
\[
 \begin{split}
 0= &\mathfrak{r}_1 E = \sum_{S_4} \operatorname{S}_{ij} \bigg( \\
  &2 c_1 \la \p^1 \p^2 \p^j \ra \la \p^3 \p^4 \p^{\mu} \ra
   \la \p^i \p^{\mu} \p^{\nu} \ra \la \p^{\nu} \ra_1 \\
  &-c_1 \la \p^1 \p^2 \p^{\mu} \ra \la \p^3 \p^4 \p^{\nu} \ra
   \la \p^i \p^{\nu} \p^{\mu} \ra \la \p^{j} \ra_1 \\
  &+c_2 \la \p^1 \p^2 \p^{\mu} \ra \la \p^3 \p^j \p^{\mu} \ra
   \la \p^i \p^4  \p^{\nu} \ra \la \p^{\nu} \ra_1 \\
  &- c_2 \la \p^1 \p^2 \p^{\mu} \ra \la \p^3 \p^{\mu} \p^{\nu} \ra
   \la \p^i \p^4  \p^{\nu} \ra \la \p^{j} \ra_1 \\
  &+ c_3 \la \p^1 \p^2 \p^{j} \ra \la \p^3 \p^4 \p^{i}_1 \p^{\nu} \ra
   \la \p^{\nu} \ra_1 \\
  &+ c_3 \la \p^1 \p^2 \p^{\mu} \ra \la \p^3 \p^4 \p^{i}_1 \p^{\nu} \ra
   \la \p^{j} \ra_1 \\
  &- c_3 \la \p^1 \p^2 \p^{\mu} \ra \la \p^3 \p^4 \p^{i} \ra
   \la \p^j \p^{\mu}  \p^{\nu} \ra \la \p^{\nu} \ra_1 \\
  &-2 c_3 \la \p^1 \p^2 \p^{\mu} \ra \la \p^3 \p^{i}  \p^{\mu} \ra
   \la \p^4 \p^j \p^{\nu} \ra \la \p^{\nu} \ra_1 \\
  &+ c_4 \la \p^1 \p^2 \p^3 \p^i_1 \ra \la \p^4 \p^{j} \p^{\nu} \ra
   \la \p^{\nu} \ra_1 \\
  &-3 c_4 \la \p^1 \p^2 \p^{i} \ra \la \p^3 \p^{j} \p^{\mu} \ra
   \la \p^4 \p^{\mu} \p^{\nu} \ra \la \p^{\nu} \ra_1 \bigg) \\
  &+\text{genus-zero-only terms} .
 \end{split}
\]
Here $\operatorname{S}_{ij}$ is the symmetrization operator of the indices 
$i j$.

\subsubsection{Step 3: genus one terms}
The basic strategy is to find a basis,  express the vector in terms of
the basis,  and set the coefficients to $0$.

It is easy to see that the terms containing $\la \p^j \ra_1$ gives the 
condition (after applying genus zero TRR)
\[
 -c_1 -c_2 +c_3 =0.
\]
The terms containing $\la \p^*  \p^{**}  \p^j \ra \la \p^{\nu} \ra_1$ gives
the equation
\[
 2 c_1 -3 c_4 =0.
\]
The terms containing $\la \p^*  \p^{**}  \p^{\nu} \ra \la \p^{\nu} \ra_1$ 
gives the equation
\[
 c_2 -2 c_3 + c_4 =0.
\]

\subsubsection{Step 4: Genus zero terms}
For the terms involving geometric genus zero graphs only, the only relations  
are WDVV, after stripping off all descendents by genus zero TRR.

(a) Those terms containing a factor $\la \p^*, \p^{**},\p^{***}, \p^i \ra$ give
the equation
\[
 \begin{split}
 \sum_{S_4} \operatorname{S}_{ij} &\la \p^1 \p^{2} \p^{3} \p^i \ra \Bigg[
  c_5 \la \p^4  \p^j_1 \p^{\nu} \p^{\nu} \ra \\
  &-4 c_6 \la \p^4  \p^j \p^{\mu} \ra \la \p^{\mu} \p^{\nu} \p^{\nu} \ra \\
  &-c_9 \la \p^4 \p^{\mu} \p^{\nu} \ra \la \p^{j} \p^{\mu} \p^{\nu} \ra 
  \Bigg]=0,
 \end{split}
\]
which gives condition
\[
 c_5 -4 c_6 -c_9 =0.
\]

(b) Those terms containing a factor $\la \p^*  \p^{**} \p^i \ra$ give the 
equation
\begin{equation} \label{e:5}
 \begin{split}
 \sum_{S_4} \operatorname{S}_{ij} \la \p^1 \p^{2} \p^i \ra \Bigg[
  &\frac{1}{12} c_1 \la \p^3 \p^4 \p^{\nu}\ra
   \la \p^j \p^{\mu} \p^{\mu} \p^{\nu} \ra \\
  &-3 c_5 \la \p^3 \p^j \p^{\nu}\ra \la \p^4 \p^{\mu} \p^{\mu} \p^{\nu} \ra \\
  &-6 c_6 \la \p^3 \p^4 \p^j \p^{\nu}\ra \la \p^{\mu} \p^{\mu} \p^{\nu} \ra \\
  &-2 c_7 \la \p^3 \p^4 \p^{\mu} \p^{\nu}\ra
   \la \p^{j} \p^{\mu} \p^{\nu} \ra \\
  &+ c_8 \la \p^3 \p^4 \p^j_1 \p^{\nu} \p^{\nu} \ra \\
  &- c_8 \la \p^3 \p^4 \p^{\mu} \ra \la \p^{j} \p^{\mu} \p^{\nu} \p^{\nu} \ra\\
  &-3 c_9 \la \p^3 \p^{\mu} \p^{\nu}\ra
   \la \p^4 \p^{j} \p^{\mu} \p^{\nu} \ra \Bigg] =0. \\
 \end{split}
\end{equation}

All graphs are disconnected with two components. 
One common connected component is $\la \p^1 \p^{2} \p^i \ra$. 
The other one is a stratum isomorphic to 
$\Mbar_{0,3} \times \Mbar_{0,4} \subset \Mbar_{0,5}$. 
Due to different labeling, there are 5 different strata for the second 
component, up to the obvious permutation symmetry in $\p^3,\p^4$.
\[
 \begin{split}
 \vec{v}_1 =\la \p^3 \p^4 \p^{\nu} \ra \la \p^j \p^{\mu} \p^{\mu} \p^{\nu} \ra,\\
 \vec{v}_2 =\la \p^3 \p^j \p^{\nu} \ra \la \p^4 \p^{\mu} \p^{\mu} \p^{\nu} \ra, \\
\vec{v}_3 =\la \p^3 \p^{\mu} \p^{\nu} \ra \la \p^4 \p^j \p^{\mu} \p^{\nu} \ra,\\
 \vec{v}_4 =\la \p^j \p^{\mu} \p^{\nu} \ra \la \p^3 \p^4 \p^{\mu} \p^{\nu} \ra,\\
\vec{v}_5 =\la \p^{\mu} \p^{\mu} \p^{\nu} \ra \la \p^3 \p^4 \p^j \p^{\nu} \ra.
 \end{split}
\]

The WDVV equations induce three linear relations in 
$\vec{v}_1, \ldots, \vec{v}_5$; two of them are independent.
\[
 \begin{split}
  &\vec{v}_1 + \vec{v}_5 = 2 \vec{v}_3, \\
  &\vec{v}_2 + \vec{v}_5 = \vec{v}_3 + \vec{v}_4, \\
  &\vec{v}_1 + \vec{v}_4 = \vec{v}_2 + \vec{v}_3.
 \end{split}
\]
Thus, one can write $\vec{v}_1$ and $\vec{v}_2$ in terms of 
$\vec{v}_3, \vec{v}_4,\vec{v}_5$.
Equation \eqref{e:5} then gives
\[
 \begin{split}
 \frac{1}{6} c_1 - 3 c_5 + 3 c_9=0, \\
 -3 c_5 -2 c_7 + 2 c_8 =0, \\
 - \frac{1}{12} c_1 + 3 c_5 - 6 c_6 =0.
 \end{split}
\]

(c) The remaining terms, after genus zero TRR,
contain no descendents. Therefore
the only relations are WDVV and their derivatives.
However, WDVV and their derivatives don't change the summation of the
coefficients, therefore the summation has to vanish. This gives another equation
\[
 -\frac{1}{2} c_1 - \frac{11}{24} c_2 -\frac{11}{24} c_3 - \frac{11}{24} c_4
 + 3 c_6 - 3 c_8 =0.
\]

\subsubsection{Step 4: Final equation}
Combining the above equations, one can express all coefficients
in terms of $c_3$ and $c_9$:
\[
 \begin{split}
  &c_1 = -3 c_3, \ c_2=4 c_3, \ c_4= -2 c_3, \ c_5= -\frac{1}{6}c_3 - c_9, \\
  &c_6 = -\frac{1}{24}c_3 - \frac{1}{2}c_9, \ c_7 = \frac{1}{4}c_3 + c_9, \
   c_8 = - \frac{1}{2}c_9.
 \end{split}
\]
That is,
\[
  E= - c_3 (\text{Getzler's coefficients}) + c_9 (T) =0,
\]
where $T$ is a sum of (geometric) genus zero graphs. 
It is easy to see that $T=0$ by WDVV.
Therefore, $l=1$ case is established.

\subsubsection{Step 6: $l=2$}

By the same computation, one is led to
\[
 \begin{split}
 \mathfrak{r_2} E = \sum_{S_4} \operatorname{A}_{ij} \Bigg[
 &6 \la \p^1 \p^2 \p^j \ra \la \p^3 \p^4 \p^{\nu} \ra \la \p^i_2 \p^{\nu} \ra_1\\
 &-4 \la \p^1 \p^2 \p_{\mu} \ra \la \p^3 \p^{\mu} \p^j \ra \la \p^4 \p^i_2 \ra_1\\
 &+\frac{1}{24} \la \p^1 \p^2 \p^3 \p^4 \p^i_2 \ra \la \p^j \p^{\nu} \p^{\nu} \ra \\
 &+3\la \p^1 \p^2 \p^{\mu} \ra \la \p^3 \p^4 \p^{\nu} \ra
  \la \p^i \p^{\mu} \p^{\nu} \ra \la \p^j_1 \ra_1 \\
 &-4\la \p^1 \p^2 \p^{\mu} \ra \la \p^3 \p^{\mu} \p^{\nu} \ra
  \la \p^i \p^4 \p^{\nu} \ra \la \p^j_1 \ra_1 \\
 &-\frac{4}{24} \la \p^i_1 \p^1 \p^2 \p^3 \ra \la \p^j \p^4 \p^{\nu} \ra
  \la \p^{\mu} \p^{\nu} \p^{\nu} \ra \\
 &-\frac{6}{24} \la \p^i_1 \p^1 \p^2 \p^{\mu} \ra \la \p^j \p^3 \p^4 \ra
  \la \p^{\mu} \p^{\nu} \p^{\nu} \ra  \Bigg] ,
 \end{split}
\]
where $\operatorname{A}_{ij}$ is the anti-symmetrizer.
By easy application of genus one TRR and WDVV and the antisymmetric
property of $i,j$, the first term cancels with the seventh term;
the second term cancels with the sixth term; the third, fourth and
fifth terms combine to vanish.

By Lemma~\ref{l:1}, it is enough to check $R$-invariance for $l=1$ and $l=2$.
The proof is complete.

\subsection{Final remarks}
The three Conjectures proposed fit very well with the current understanding
of the tautological rings. 
For example, they are consistent with Graber--Vakil's stratification of the 
tautological rings by the number of rational components and with the 
Poincar\'e duality conjecture. 
See \cite{rV} for the statements and references.

As for the proof of the conjectures, there are possibly a few
approaches. One approach is to go through its connection to the
Gromov--Witten theory. The key idea is to utilize some distinguished
points on the ``moduli spaces'' of axiomatic GW theories, for which
the theories and the ``quantized loop group'' elements moving them
around are well studied. For example, let $\tau^{pt}$ be the GW
theory of $N$ points and $\tau^{\PP^{N-1}}$ be the GW theory of
$\PP^{N-1}$. There is a ``quantized loop group'' element $\hat{R}_N$
such that
\[
 \hat{R}_N \tau^{pt} = \tau^{\PP^{N-1}}.
\]
Since the tautological equations certainly hold for both theories,
they are invariant under the action of $\hat{R}_N$.
It is possible that the invariance under enough group elements like $\hat{R}_N$
implies invariance under the whole ``quantized loop group''.

Another approach is to consider the invariance condition as some hidden
geometric structure on the geometry of the moduli spaces of curves.
One might be able to unravel this structure.

\end{document}